\documentclass[11pt,leqno,oneside,letterpaper]{amsart}

\usepackage[T1]{fontenc}
\usepackage[utf8]{inputenc}
\usepackage{lmodern}
\usepackage{eulervm}
\usepackage{amsmath,amssymb,amsthm,mathtools}
\usepackage[margin=1.15in]{geometry}
\usepackage{booktabs,array,longtable}
\usepackage{enumitem}
\usepackage{xurl}
\usepackage[colorlinks,linkcolor=blue,citecolor=blue,urlcolor=blue]{hyperref}
\usepackage{orcidlink}
\setlist{leftmargin=*,itemsep=3pt,topsep=5pt}
\newtheorem{theorem}{Theorem}[section]
\newtheorem{proposition}[theorem]{Proposition}
\newtheorem{lemma}[theorem]{Lemma}
\newtheorem{corollary}[theorem]{Corollary}

\theoremstyle{definition}

\theoremstyle{remark}

\numberwithin{equation}{section}
\newcommand{\Z}{\mathbb Z}
\newcommand{\Q}{\mathbb Q}
\newcommand{\RR}{\mathbb R}
\newcommand{\OK}{\mathcal O_K}
\newcommand{\Nm}{\operatorname{N}}
\newcommand{\Tr}{\operatorname{Tr}}
\newcommand{\disc}{\operatorname{disc}}

\newcommand{\gp}[1]{\langle#1\rangle}

\title[Affine pairs of minimizing units]{Affine pairs of minimizing units in totally real non-Galois cubic fields}
\author[J.~Lu]{Junyu Lu\,\orcidlink{0009-0001-4973-3198}}
\address{Department of Mathematics, Yang-En University, Quanzhou, Fujian, China, 362014}
\email{jy.lu@outlook.com}
\date{}
\subjclass[2020]{11R16, 11R27, 11H06}
\keywords{totally real cubic field, minimizing unit, unit group, affine equivalence, exceptional unit}
\hypersetup{pdftitle={Affine pairs of minimizing units in totally real non-Galois cubic fields},pdfauthor={Junyu Lu}}

\begin{document}
\begin{abstract}
  Let $\Gamma$ be a rank-two subgroup of the integral norm-one units in a totally real non-Galois cubic field. We prove that two distinct affinely equivalent minima of Godwin's spread on $\Gamma\setminus\{1\}$ generate $\Gamma$ and determine its complete minimum set: a pair with sum $\pm2$ or a triple of integral translates of product one. Every minimum generates the order $\Z[\Gamma]$. For an order whose norm-one units have such a pair, the order generated by its units is the smallest suborder with the same minimum value. We also give explicit families and a maximal order with minima in two distinct affine classes.
\end{abstract}
\maketitle

\section{Introduction and main results}\label{sec:introduction}

Let $K$ be a totally real cubic field with ordered real embeddings $\sigma_1,\sigma_2,\sigma_3:K\hookrightarrow\RR$, and write $x_i=\sigma_i(x)$. For an order $\Omega\subset K$, put
\[
  E^1(\Omega)=\{u\in\Omega^\times:\Nm_{K/\Q}(u)=1\},\qquad
  S(x)=\frac12\sum_{i<j}(x_i-x_j)^2
  =\frac{3\Tr_{K/\Q}(x^2)-\Tr_{K/\Q}(x)^2}{2}.
\]
The formula $S(x_1,x_2,x_3)=\frac12\sum_{i<j}(x_i-x_j)^2$ also defines $S$ on arbitrary real triples and the function $S$ is Godwin's spread \cite[\S2]{Godwin1960}, also studied by Gras \cite[\S1]{Gras1980} and many others.

For a rank-two subgroup $\Gamma\subseteq E^1(\Omega)$, define
\begin{equation}\label{eq:minimum-definition}
  m_\Gamma=\min_{u\in\Gamma\setminus\{1\}}S(u),\qquad
  \mathcal M_\Gamma=\{u\in\Gamma\setminus\{1\}:S(u)=m_\Gamma\}.
\end{equation}
When $\Gamma=E^1(\Omega)$, we write $m_\Omega$ and $\mathcal M_\Omega$ as well. The minimum exists and its minimizing set is finite, as shown in Section~\ref{sec:geometry}. Two algebraic integers $u,v\in K$ are called \emph{affinely equivalent} if $v=\pm u+n$ for some $n\in\Z$. Affine equivalence preserves $S$.

The map $\varepsilon\{\pm1\}\mapsto\Nm_{K/\Q}(\varepsilon)\varepsilon$ identifies $\Omega^\times/\{\pm1\}$ with $E^1(\Omega)$ and preserves $S$. Thus norm normalization removes the sign ambiguity in minimizing over units other than $\pm1$. We write $\Z[\Gamma]$ for the order generated by the units in $\Gamma$.

\begin{theorem}\label{thm:subgroup-affine-minima}
  Suppose $K$ is a totally real non-Galois cubic field and $\Gamma\subseteq E^1(\OK)$ is a rank-two subgroup. If distinct $u,v\in\mathcal M_\Gamma$ are affinely equivalent, then $\Gamma=\gp{u,v}$ and exactly one of the following holds:
  \begin{enumerate}[label=\textup{(\roman*)}]
    \item $\mathcal M_\Gamma=\{u,v\}$ and $u+v=\pm2$;
    \item $\mathcal M_\Gamma$ is a triple of distinct integral translates of product one, and every two of them generate $\Gamma$.
  \end{enumerate}
  Moreover, $\Z[\Gamma]=\Z[w]$ for every $w\in\mathcal M_\Gamma$.
\end{theorem}

An order of the form $\Z[\alpha]$, with $\alpha$ a nonrational algebraic integer in $K$, is called an \emph{equation order}. For the full norm-one unit group of an order, Theorem~\ref{thm:subgroup-affine-minima} identifies exactly which suborders preserve the minimum.

\begin{corollary}\label{cor:unit-generated-order}
  Let $\Omega$ be an order in a totally real non-Galois cubic field. Suppose two distinct elements of $\mathcal M_\Omega$ are affinely equivalent, and put $R_U=\Z[\Omega^\times]=\Z[E^1(\Omega)]$. Then $R_U=\Z[w]$ for every $w\in\mathcal M_\Omega$, and the class of $w$ is the unique shortest nonzero class modulo sign in $R_U/\Z$, with respect to $S$. A suborder $\Xi\subseteq\Omega$ satisfies $m_\Xi=m_\Omega$ if and only if $E^1(\Xi)=E^1(\Omega)$; these conditions are also equivalent to $R_U\subseteq\Xi$. Whenever they hold, $\mathcal M_\Xi=\mathcal M_\Omega$.
\end{corollary}

For any order $\Omega$, the ring $\Z[\Omega^\times]$ is the smallest suborder with the same unit group, since it contains every unit and its inverse. Under the affine-pair hypothesis, Corollary~\ref{cor:unit-generated-order} adds that every minimum generates this ring and that equality of minimum values forces equality of unit groups. The corollary is proved in Section~\ref{sec:applications}.

In a totally real non-Galois cubic field, if $\mathcal M_{\OK}$ contains several affine classes, each is a singleton by Theorem~\ref{thm:subgroup-affine-minima}: a second minimum in any one class would force all minima into that class. This possibility occurs in the field defined by $X^3-14X^2-14X+14$, whose maximal order has two inequivalent minima at spread $952$ (Section~\ref{sec:counterexample}). Thus equal spread alone does not imply affine equivalence.

For the polynomial families in Section~\ref{sec:geometry}, we show that the displayed minima remain globally minimizing in an intermediate unit group precisely when they generate that group (Theorem~\ref{thm:family-generation}).

\subsection{Relation to earlier work}
Godwin introduced the spread $S$ in his study of cubic units \cite[\S2]{Godwin1960}. Brunotte and Halter-Koch bound the index of successive $S$-minima in the full maximal-order unit group by four \cite[Satz~1]{BHK1979}. Gras proved generation by a minimum and a distinct Galois conjugate in cyclic cubic fields of conductor different from $7$ and $9$ \cite[Theorem~2]{Gras1980}. Ennola proves generation by successive $S$-minima for sufficiently large field discriminant and expresses a first minimum in terms of Berwick units \cite[Theorems~1 and~2]{Ennola1987}. Subgroup index methods also appear in \cite[Lemma~3 and Satz~2]{BHK1981}, for functions of the absolute values of conjugates.

The polynomial families and their standard units are classical. For the family~\eqref{eq:family}, Ennola records the standard triple and its spread \cite[Lemma~6(ii)]{Ennola1987}. Watabe proves minimality and generation for parameters $a\ge5$ with squarefree polynomial discriminant \cite[Theorem and \S2]{Watabe1983III}, and generation for $a=3,4$ \cite[\S3(ii)]{Watabe1984VI}. Minemura gives the affine normal forms and quadratic-form arguments under arithmetic maximality hypotheses \cite[Proposition~1, Theorems~1, 3--4 and Lemmas~5--6]{Minemura1998}. Stender proves generation for product families under separation hypotheses \cite[\S2, Satz~7]{Stender1972}, and Thomas treats the orders defined by~\eqref{eq:family} and broader product families \cite[Proposition~(3.6) and Theorem~(3.9)]{thomas1979}. We reconstruct the equation-order minimum calculations with their complete equality cases.

Komatsu classifies pairs of cubic units with nonzero rational difference in cyclic cubic fields, without a spread-minimality hypothesis \cite[Theorems~1.1 and~1.2]{komatsu}. Theorem~\ref{thm:subgroup-affine-minima} concerns global minima in prescribed subgroups of non-Galois cubic fields.

Proper inclusions of cubic orders with the same unit group already follow from Minemura's work \cite[Example~1]{minemura}. He constructs totally real cubic equation orders $\Z[\alpha]$ whose norm-one units are generated by $7\alpha+11$ and $7\alpha+43$. These generators and their inverses lie in $\Z[7\alpha]$, so $\Z[7\alpha]\subset\Z[\alpha]$ has index $343$ and the same unit group. David and Shapira construct arbitrarily long finite chains of proper inclusions with the same unit group \cite[discussion preceding Theorem~3.14 and its proof]{DavidShapira2018}. Theorem~\ref{thm:index-two-family} gives an explicit family of index two, with its conductor, a maximality criterion, and the complete minimum set.


\subsection{Organization of the article} Section~\ref{sec:geometry} gives the affine normal forms and equation-order minima. Section~\ref{sec:affine-saturation} proves the result on subgroups, including the small cases required by the relative index bound. Section~\ref{sec:applications} derives the order and family consequences and gives the counterexample.

\section{Affine pairs and equation-order minima}\label{sec:geometry}

\subsection{The quadratic form on the additive quotient}

Since $S$ is unchanged by integral translation, it is a quadratic form on the additive quotient $\Omega/\Z$. Shortest classes in this quotient give a lower bound for unit spread. Solving $\Nm_{K/\Q}(x)=1$ in those classes then determines whether the bound is attained by norm-one units and identifies every equality case.

Let $\iota:K\longrightarrow\RR^3$, $\iota(x)=(\sigma_1(x),\sigma_2(x),\sigma_3(x))$, be the Minkowski embedding and $\mathbf{1}=(1,1,1)$. Set
\[
  V=\{x\in\RR^3:x_1+x_2+x_3=0\},\qquad
  \pi(x)=x-\tfrac13(x_1+x_2+x_3)\mathbf{1}.
\]
The projected additive lattice is $\Lambda_\Omega^{\mathrm{add}}=\pi(\iota(\Omega))$. It is naturally isomorphic to $\Omega/\Z$, since $\Omega\cap\Q=\Z$. For $w\in\Omega$, write $[w]=w+\Z$ for its additive class. If $\{1,\omega_1,\omega_2\}$ is an integral basis of $\Omega$, the projections of $\omega_1,\omega_2$ are linearly independent in $V$ and generate $\Lambda_\Omega^{\mathrm{add}}$; independence follows from the nonzero determinant of the embedding matrix of $1,\omega_1,\omega_2$. The projected lattice is not the integral trace-zero submodule of $\Omega$. Godwin's search already uses the translation invariance of $S$ and its positive definite binary form on $\Omega/\Z$ \cite[\S4]{Godwin1960}. The Euclidean metric on $V$ satisfies
\begin{equation}\label{eq:projected-norm}
  S(\alpha)=\frac{3\Tr_{K/\Q}(\alpha^2)-\Tr_{K/\Q}(\alpha)^2}{2}
  =\frac32\|\pi(\iota(\alpha))\|^2.
\end{equation}
For $\alpha\in\Omega$, put $t=\Tr_{K/\Q}(\alpha)$ and $s=\sum_{i<j}\sigma_i(\alpha)\sigma_j(\alpha)$. Then $S(\alpha)=t^2-3s\in\Z$, so $S$ is an integral positive definite binary quadratic form on $\Omega/\Z$. The only rational norm-one unit is $1$, so $S$ takes positive integer values on $\Gamma\setminus\{1\}$ and the minimum in~\eqref{eq:minimum-definition} exists. A bounded disk contains finitely many projected lattice points. In each integral-translation class, the equation $\Nm_{K/\Q}(\alpha+n)=1$ is cubic in $n$ and has at most three solutions. Thus only finitely many norm-one units have bounded spread, and $\mathcal M_\Gamma$ is finite.

At a fixed spread, affine equivalence can also be recognized from the generated orders. Here the non-Galois hypothesis is essential.

\begin{proposition}\label{prop:equivalent-orders}
  Let $\Omega$ be an order in a totally real non-Galois cubic field $K$. Suppose $\alpha,\beta\in\Omega\setminus\Z$ and $S(\alpha)=S(\beta)$. The following are equivalent:
  \begin{enumerate}[label=\textup{(\roman*)}]
    \item $\beta=\pm\alpha+n$ for an integer $n$;
    \item $\Z[\alpha]=\Z[\beta]$ as embedded subrings of $K$;
    \item $[\Omega:\Z[\alpha]]=[\Omega:\Z[\beta]]$.
  \end{enumerate}
\end{proposition}
\begin{proof}
  The first assertion implies the second, and the second implies the third. Equal indices imply equal polynomial discriminants, by the change-of-basis formula for the trace determinant. If $s=S(\alpha)=S(\beta)$, the centered elements $\alpha-\Tr_{K/\Q}(\alpha)/3$ and $\beta-\Tr_{K/\Q}(\beta)/3$ have minimal polynomials $X^3-(s/3)X+b_\alpha$ and $X^3-(s/3)X+b_\beta$. The discriminants are $4s^3/27-27b_\alpha^2$ and $4s^3/27-27b_\beta^2$. Hence $b_\alpha=\pm b_\beta$. After a possible negation, the centered elements are roots of the same irreducible cubic. A non-Galois cubic field cannot contain two distinct roots of such a cubic: their sum gives the third, which would make the field normal. The centered elements therefore agree. Their difference before centering is a rational algebraic integer, hence an integer.
\end{proof}

\subsection{The Ennola family}
The Ennola family of polynomials gives a first example in which the shortest additive classes contain units. The quadratic form induced by $S$ also gives the complete minimum set without assuming that the equation order is maximal.
For an integer $a\ge3$, the \emph{Ennola polynomial} with parameter $a$ is
\begin{equation}\label{eq:family}
  f_a(X)=X^3+(a-1)X^2-aX-1=X(X-1)(X+a)-1.
\end{equation}
Let $\theta$ be its greatest real root and put $K_a=\Q(\theta)$, $R_a=\Z[\theta]$, $\mathcal O_a=\mathcal O_{K_a}$ and $H_a=\gp{\theta,\theta-1}$. We refer to $K_a$ as an \emph{Ennola field}; $R_a$ is the equation order defined by $f_a$. Then
\begin{equation}\label{eq:exceptional-product}
  \theta(\theta-1)(\theta+a)=1,
\end{equation}
and all three factors have norm one. In particular, $\theta$ is an exceptional unit, that is, both $\theta$ and $\theta-1$ are units.
The rational-root test gives irreducibility of $f_a$, since $f_a(1)=-1$ and $f_a(-1)=2a-3$. Its signs at $-a,1-a,-1,0,1,2$ show that its three roots lie respectively in
\[
  (-a,1-a),\qquad(-1,0),\qquad(1,2).
\]
Expansion gives
\begin{equation}\label{eq:discriminants}
  \Delta_a:=\disc(f_a)=(a^2+3a-1)^2-32
    =[\mathcal O_a:R_a]^2D_{K_a}.
\end{equation}
For $a\ge3$, one has $a^2+3a-1\ge17$, and $(a^2+3a-2)^2<\Delta_a<(a^2+3a-1)^2$. Thus $\Delta_a$ is not a square and $K_a$ is non-Galois. In~\eqref{eq:discriminants}, $\Delta_a$ is the polynomial discriminant, also equal to $\disc(R_a)$, whereas $D_{K_a}=\disc(\mathcal O_a)$ is the field discriminant.

Thomas's equation-order generation result uses $\lambda=1-1/\theta$ and $n=a+4$ \cite[Proposition~(3.6)]{thomas1979}, as in \cite{Louboutin2017}. Since $1-\lambda$ is a unit and $\theta=(1-\lambda)^{-1}$, one has $\Z[\lambda]=R_a$. Thomas's generators $\lambda,\lambda-1$ have norm-one normalizations $(\theta-1)/\theta,\theta^{-1}$, which generate $H_a$. In Minemura's notation, the generator is $\theta-1$ and the parameters are $t=1$, $s=a+1$ \cite{Minemura1998}.

\begin{theorem}\label{thm:equation-order}
  For every $a\ge3$, the minimum of $S$ on $E^1(R_a)\setminus\{1\}$ is $a^2+a+1$, and its complete set of minimizing units is
  \[
    \{\theta-1,\theta,\theta+a\}.
  \]
  The same assertion holds in every group $H_a\subseteq\Gamma\subseteq E^1(R_a)$.
\end{theorem}

\begin{proof}
  For $a\ge5$, Watabe's inequalities in the basis $1,\delta,\delta^2$, with $\delta=\theta-1$ and $n=-a-2$, already show that $\pm\theta+\Z$ are the only classes of least positive spread in $R_a/\Z$ \cite[\S2]{Watabe1983III}. We use a different basis to handle all $a\ge3$ in the same calculation.

  Put $\eta=\theta^2+(a-1)\theta=\theta^{-1}+a$. Then $1,\theta,\eta$ is an integral basis of the equation order, and $\theta\eta=1+a\theta$. Newton's identities give
  \begin{gather*}
    \Tr_{K_a/\Q}(\theta)=1-a,\quad \Tr_{K_a/\Q}(\theta^2)=a^2+1,\\
    \Tr_{K_a/\Q}(\eta)=2a,\quad \Tr_{K_a/\Q}(\eta^2)=2a^2+2a-2,\quad
    \Tr_{K_a/\Q}(\theta\eta)=3+a-a^2.
  \end{gather*}
  Consequently
  \begin{equation}\label{eq:equation-form}
    S(b\theta+c\eta)=Ab^2+Bbc+Cc^2,
    \quad A=a^2+a+1,\quad B=-a^2+a+9,\quad C=a^2+3a-3.
  \end{equation}
  For $a\ge3$, $C>A>|B|$. If $b,c\ne0$, put $x=|b|$, $y=|c|$, $d=|B|$. Then
  \[
    S(b\theta+c\eta)\ge(A-d/2)x^2+(C-d/2)y^2+(d/2)(x-y)^2
    \ge A+C-d>A.
  \]
  If $b=0,c\ne0$, the value is greater than $A$; if $c=0$, equality with $A$ holds only for $b=\pm1$. The elements of least positive $S$ are therefore exactly $\pm\theta+n$. Their norm equations are
  \[
    \Nm_{K_a/\Q}(\theta+n)-1=n(n+1)(n-a),\qquad
    \Nm_{K_a/\Q}(-\theta+n)-1=n(n-1)(n+a)-2.
  \]
  The first gives $n=-1,0,a$. For the second, $n=0,1,-a$ gives product zero; $n\ge2$ gives product at least ten; $n\le-a-1$ gives a negative product; and $n=-k$, $1\le k\le a-1$, gives $k(k+1)(a-k)$, which could equal two only for $k=1,a=2$. Thus the second has no integer solution. Identity~\eqref{eq:exceptional-product} places all three remaining units in $H_a$.
\end{proof}

If $R_a=\mathcal O_a$, the theorem gives the complete field minimum; squarefree $\Delta_a$ ensures this equality by~\eqref{eq:discriminants}.

\subsection{Powers and independence of minima}\label{sec:eventual-powers}

Throughout this subsection, $K$ is a totally real non-Galois cubic field. To apply the index bounds in Section~\ref{sec:affine-saturation}, we need distinct minima to be multiplicatively independent. We prove this by showing that taking a proper positive power increases spread and that a unit and its inverse have different spreads. A nonidentity element of a free abelian group is called primitive if it is not a proper power in that group.

\begin{lemma}\label{lem:eventual-small-S}
  If $\alpha$ is a nonrational algebraic integer in a non-Galois totally real cubic field, then $S(\alpha)\ge10$.
\end{lemma}
\begin{proof}
  Translate $\alpha$ by an integer so that its trace belongs to $\{-1,0,1\}$, and negate it if that trace is $-1$. These operations preserve $S$ and the generated field. The resulting integer has irreducible polynomial
  \[
    g(X)=X^3-tX^2+sX+c,\qquad t\in\{0,1\},\quad s,c\in\Z,
  \]
  with
  \[
    0<S=t^2-3s,\qquad \disc(g)=t^2s^2-4s^3+(4t^3-18ts)c-27c^2>0.
  \]
  If $S\le9$, then $s\in\{-3,-2,-1\}$ for $t=0$ and $s\in\{-2,-1,0\}$ for $t=1$. Positivity of the displayed concave quadratic in $c$ forces $|c|\le2$. Checking these integers leaves only the irreducible polynomials $X^3-3X\pm1$ and $X^3-X^2-2X+1$, of discriminants $81$ and $49$.

  Every irreducible possibility therefore has square discriminant. To see the contradiction directly, if $\beta,\beta',\beta''$ are its roots and its discriminant is $d^2$ with $d\in\Q$, then $\beta'-\beta''=\pm d/g'(\beta)$ and $\beta'+\beta''=t-\beta$ belong to $\Q(\beta)$. Thus all roots belong to the cubic field, making it Galois. Hence $S>9$, and integrality gives the assertion.
\end{proof}

Let $w\in E^1(\OK)\setminus\{1\}$. Godwin's square estimate for integral cubic units \cite[\S3, discussion following Lemma~3]{Godwin1960} is
\[
  S(w)^2\le9S(w^2).
\]
Lemma~\ref{lem:eventual-small-S} therefore gives $S(w^2)\ge S(w)^2/9>S(w)$. For odd $k\ge3$ and real $r,s$,
\[
  |r^k-s^k|\ge2^{1-k}|r-s|^k.
\]
For equal signs the binomial theorem gives the stronger coefficient one; for opposite signs the bound follows from convexity on the nonnegative reals. Applying the displayed inequality to the three conjugate pairs and then using convexity gives
\[
  S(w^k)\ge\frac{2^{2-2k}}2\,3^{1-k}
  \left(\sum_{i<j}|w_i-w_j|^2\right)^k
    =6^{1-k}S(w)^k>S(w).
\]
Write $k=2^r\ell$ with $\ell$ odd. The odd-power estimate and successive squaring give $S(w^k)>S(w)$ for every $k\ge2$; every intermediate power is again a nonidentity norm-one unit. Consequently, a minimizing unit is primitive in the group over which the minimum is taken. To compare $w$ with its inverse, put $t=\Tr_{K/\Q}(w)$ and $s=\Tr_{K/\Q}(w^{-1})$. Then
\begin{equation}\label{eq:eventual-inverse-difference}
  S(w)-S(w^{-1})=(t-s)(t+s+3)\ne0.
\end{equation}
Indeed, if $t=s$, then $1$ is a root of $X^3-tX^2+tX-1$, contradicting irreducibility. If $t+s+3=0$, then $\disc(X^3-tX^2+sX-1)=(t^2+3t+9)^2$, contradicting the non-Galois hypothesis. In a free abelian group, two primitive elements on the same rational line are equal or inverse. Thus two distinct minimizing units in a rank-two group are multiplicatively independent. Dirichlet's unit theorem gives rank two for the norm-one unit group of a totally real cubic field; its torsion is trivial because the only real roots of unity are $\pm1$ and $\Nm_{K/\Q}(-1)=-1$.

The norm-one unit group of every order $\Omega$ has finite index in $E^1(\OK)$ and hence rank two. Indeed, for a nonzero ideal $\mathfrak f$ of $\OK$ contained in $\Omega$, the kernel of reduction $E^1(\OK)\to(\OK/\mathfrak f)^\times$ lies in $E^1(\Omega)$: both a unit in this kernel and its inverse belong to $1+\mathfrak f\subseteq\Omega$. The target is finite.

\subsection{Affine unit classes}\label{sec:translates}

Throughout this subsection $K$ is a totally real non-Galois cubic field. We classify the norm-one units in a single affine class before imposing minimality. For a nonrational algebraic integer $\beta$, put $F(X)=\Nm_{K/\Q}(\beta+X)$. The integer solutions of $F(n)=\pm1$ correspond bijectively to norm-one units in the affine class $\{\pm\beta+n:n\in\Z\}$, via $n\mapsto F(n)(\beta+n)$. Every norm-one unit in the class has this form. Two solutions with the same norm sign give distinct units, while equality across opposite signs would imply $2\beta\in\Z$.

\begin{lemma}\label{lem:translate-triple}
  An affine class contains at most three norm-one units. If it contains three, they are distinct integral translates of product one. Writing them as $\theta-b,\theta,\theta+a$, with $a,b$ positive integers, one has
  \[
    \theta(\theta+a)(\theta-b)=1,\qquad S(\theta)=a^2+ab+b^2.
  \]
\end{lemma}
\begin{proof}
  If $F(r)=1$ and $F(s)=-1$, then $r-s\mid2$. Suppose three distinct unit translates have mixed norm signs. Translate the singleton sign to zero, obtaining $F(0)=-\varepsilon$ and $F(r)=F(s)=\varepsilon$, where $\varepsilon=\pm1$ and $r,s$ are distinct elements of $\{-2,-1,1,2\}$. The third root $d$ of the monic polynomial $F-\varepsilon$ is an integer, and evaluation at zero gives $rsd=2\varepsilon$. Substitution of these possibilities gives exactly
  \[
    \begin{gathered}
      X^3+4X^2+5X+1,\quad X^3-4X^2+5X-1,\quad X^3-3X\pm1,\\
      X^3+2X^2-X-1,\quad X^3-2X^2-X+1.
    \end{gathered}
  \]
  Their discriminants are respectively $-23,-23,81,81,49,49$. None can define a totally real non-Galois cubic field. Therefore every three unit translates have a common norm. Four are impossible: either a mixed triple occurs or a cubic takes one fixed value at four distinct arguments.

  If the three translates $\beta+n_i$ have common norm $\varepsilon$, then $F(X)-\varepsilon=\prod_i(X-n_i)$. Evaluating at $X=-\beta$ gives $\prod_i(\beta+n_i)=\varepsilon$. Multiplication of each factor by $\varepsilon$ normalizes its norm to one and makes the product one. The normalized elements remain literal translates. Taking the middle one as $\theta$ gives its polynomial $X(X+a)(X-b)-1$. Its trace and second symmetric function are $b-a$ and $-ab$, so its spread is $(b-a)^2+3ab=a^2+ab+b^2$.
\end{proof}

\begin{lemma}\label{lem:translate-double}
  Let $\Gamma\subseteq E^1(\OK)$ be a rank-two subgroup. If an affine class in $\mathcal M_\Gamma$ contains exactly two representatives $u,v$, then $u+v=\pm2$. There is an algebraic integer $\beta$ with $\Z[\beta]=\Z[u]$ and $T=\Tr_{K/\Q}(\beta)>0$ such that
  \begin{equation}\label{eq:midpoint-family}
    \beta^3-T\beta^2+(\delta-1)\beta+T=0,\qquad
    \{u,v\}=\{\delta(\beta+1),-\delta(\beta-1)\},
  \end{equation}
  where $\delta\in\{1,-1\}$. The common minimum is $T^2+3(1-\delta)$.
\end{lemma}
\begin{proof}
  Suppose first that $u,v$ are integral translates. The polynomial $\Nm_{K/\Q}(u+X)-1$ has two distinct integer roots, so its third root is also an integer. Evaluating its factorization at $X=-u$ shows that the third norm-one translate is $(uv)^{-1}\in\Gamma$. This translate has the same spread and hence must equal $u$ or $v$. If it equals $w$, the other unit is $w^{-2}$, contradicting the primitivity of a $\Gamma$-minimum established in Section~\ref{sec:eventual-powers}. Thus $u,v$ form a reflected pair: $v=-u+n$ for some $n\in\Z$. The two translates $u,-v$ then have opposite norms, so their translation distance divides two.

  To exclude distance one, normalize one endpoint $u$ to norm one and write $t=\Tr_{K/\Q}(u)$, $s=\Tr_{K/\Q}(u^{-1})$. If $\Nm_{K/\Q}(u+1)=-1$, then $s=-t-3$, giving square polynomial discriminant $(t^2+3t+9)^2$. If instead $\Nm_{K/\Q}(u-1)=-1$, then $s=t+1$ and the other normalized unit is $v=1-u$. In this case
  \[
    S(u)=S(v)=t^2-3t-3,\quad
    S(u^{-1})-S(u)=2t+4,\quad
    S(v^{-1})-S(v)=10-2t.
  \]
  Minimality forces $-2\le t\le5$. At these eight integers the polynomial discriminant is, in order,
  \[
    49,\ -23,\ -31,\ -23,\ -23,\ -31,\ -23,\ 49.
  \]
  All cases are excluded. The distance is therefore two.

  Let $\beta$ be the midpoint of the two translates and put $\delta=\Nm_{K/\Q}(\beta+1)=-\Nm_{K/\Q}(\beta-1)$. For this midpoint, write $F(X)=\Nm_{K/\Q}(\beta+X)=X^3+TX^2+EX+C$. The values at $1,-1$ give $E=\delta-1$ and $C=-T$. This proves the polynomial in~\eqref{eq:midpoint-family}. The trace $T$ cannot vanish, for then this polynomial is reducible; replacing $\beta$ by $-\beta$ preserves $\delta$ and makes $T>0$. Norm normalization gives the stated pair. Since $\beta=\pm u\pm1$, its equation order equals $\Z[u]$.
\end{proof}

The admissible non-Galois totally real parameters in~\eqref{eq:midpoint-family} are
\begin{equation}\label{eq:midpoint-parameters}
  \delta=1,\quad T\ge4,\ T\ne7;\qquad
  \delta=-1,\quad T\ge2.
\end{equation}
We call the cases $\delta=1$ and $\delta=-1$ the plus-midpoint and minus-midpoint families, respectively.

Indeed, the discriminants are
\begin{equation}\label{eq:midpoint-discriminants}
  \Delta_{T,1}=T^2(4T^2-27),\qquad
  \Delta_{T,-1}=4T^4+13T^2+32.
\end{equation}
In the first case positivity requires $T\ge3$, and $4T^2-27$ is a square only for $T=3,7$, by factoring $(2T-h)(2T+h)=27$. In the second case $T=1$ gives $49$; for $T\ge3$ the discriminant lies strictly between $(2T^2+3)^2$ and $(2T^2+4)^2$, while $T=2$ gives $148$. Both polynomials are irreducible: an integer root $k$ with $|k|\ge2$ would require $T=k+k/(k^2-1)$ or $T=k-k/(k^2-1)$, respectively, impossible for an integer $T$. The values $k=0,\pm1$ are excluded directly.

\subsection{The remaining equation-order families}\label{sec:translate-orders}

The affine classification gives two families of equation orders. The product family consists of $\Z[\theta]$ with $\theta(\theta+a)(\theta-b)=1$, as in Lemma~\ref{lem:translate-triple}; the midpoint families consist of $\Z[\beta]$ with $\beta$ satisfying~\eqref{eq:midpoint-family}. We now determine their shortest additive classes and the norm-one units in those classes. These equality cases will determine the full minimum set once generation is proved. The quadratic-form arguments follow \cite[Lemma~5 and proofs of Theorems~3--4]{Minemura1998}. For the product family below, Thomas's generation theorem \cite[Theorem~(3.9)]{thomas1979} applies with $\lambda=\theta+a$, $r=a$ and $s=a+b$.

\begin{proposition}\label{prop:product-order}
  Suppose $a,b\in\Z$ with $a,b\ge2$ and $\theta$ satisfies $X(X+a)(X-b)-1=0$. Then in $\Z[\theta]$, the unique shortest nonzero additive class modulo sign is $[\theta]$. Its spread is $M=a^2+ab+b^2$, and its norm-one representatives are exactly $\theta-b,\theta,\theta+a$.
\end{proposition}
\begin{proof}
  The polynomial has no rational root, since its only possible integer roots are $\pm1$ and neither works. Put $\eta=\theta^2+(a-b)\theta$. Newton sums give
  \[
    S(x\theta+y\eta)=Mx^2+Bxy+Cy^2,\qquad
    B=ab(b-a)+9,\quad C=a^2b^2+3(a-b),
  \]
  with
  \[
    4MC-B^2=3\Delta_{a,b},\qquad
    \Delta_{a,b}=a^2b^2(a+b)^2+4(a-b)^3+18ab(a-b)-27.
  \]
  Put $k=\min(a,b)\ge2$ and $h=|a-b|$. The smaller discriminant occurs at $a=k,b=k+h$, so expansion gives
  \begin{equation}\label{eq:product-discriminant-margin}
    \begin{aligned}
      3\Delta_{a,b}-k^2M^2\ge{} & 3k^6-81+18k^2(k^3-3)h                  \\
                                & +6k(4k^3-9)h^2+12(k^3-1)h^3+2k^2h^4>0.
    \end{aligned}
  \end{equation}
  Viewed as a polynomial in $h$, the right-hand side has positive coefficients for $k\ge2$; its constant coefficient satisfies $3k^6-81\ge111>0$. Thus $3\Delta_{a,b}>4M^2$, the field is totally real, and
  \[
    4M S(x\theta+y\eta)=(2Mx+By)^2+3\Delta_{a,b}y^2
  \]
  makes every vector with $y\ne0$ longer than $M$. For $y=0$, equality occurs exactly at $x=\pm1$.

  Finally $\Nm_{\Q(\theta)/\Q}(\theta+n)=1+n(n-a)(n+b)$. Norm one gives $n=-b,0,a$. Norm minus one would require $n(n-a)(n+b)=-2$, hence $n\in\{-2,-1,1,2\}$. For the two negative choices the product is nonnegative; at $1$ it is at most $-3$, and at $2$ it is zero or at most $-8$. Thus no such translate exists, and negation gives no additional norm-one representative.
\end{proof}

For the minus-midpoint family, Minemura uses the generator $-\beta-1$ and parameter $t=T+1$; for the plus-midpoint family, he uses $\beta+1$ and $t=-T-1$ \cite[Theorems~3--4]{Minemura1998}.

\begin{proposition}\label{prop:midpoint-order}
  Under~\eqref{eq:midpoint-parameters}, the unique shortest nonzero additive class of $\Z[\beta]$ modulo sign is $[\beta]$. Its spread is $T^2+3(1-\delta)$, and the complete set of norm-one minimizing units is the pair in~\eqref{eq:midpoint-family}.
\end{proposition}
\begin{proof}
  Put $\eta=\beta^2-T\beta$. On the quotient basis $\beta,\eta$, Newton sums give
  \[
    S(x\beta+y\eta)=Ax^2+Bxy+Cy^2,\quad
    A=T^2+3(1-\delta),\quad B=-T(\delta+8),\quad
    C=3T^2+(\delta-1)^2.
  \]
  If $\delta=1$, then $C-B^2/(4A)=3T^2-81/4>A$ for $T\ge4$. If $\delta=-1$, the inequality $C-B^2/(4A)>A$ is equivalent to $8T^4-9T^2-48>0$, valid at $T=2$ and increasing thereafter. Completing the square in $Ax^2+Bxy+Cy^2$ shows that $y\ne0$ gives spread greater than $A$. When $y=0$, equality with $A$ holds exactly for $x=\pm1$. The translates $\beta\pm1$ have opposite unit norms, so Lemma~\ref{lem:translate-triple} excludes a third unit translate. Norm normalization gives exactly the stated pair.
\end{proof}

These calculations determine the minima in the equation orders. A larger order may contain units with smaller spread, and a larger subgroup may contain units outside the equation order. To apply the calculations to $\Gamma$, we must first prove that the given affine pair generates $\Gamma$.

\section{Proof of the subgroup theorem}\label{sec:affine-saturation}

The generation proof has two parts. Above spread $41$, a relative form of the Brunotte--Halter-Koch index argument reduces the problem to excluding indices two, three and four inside $\Gamma$. At spread at most $41$, the affine normal forms leave nineteen cases. We first prove an independent maximal-order generation lemma and construct an index-two overorder, both needed to settle those small cases. The cutoff $41$ comes from the bound $m_\Gamma<42$ in the excluded index cases and the integrality of spread.

\subsection{The maximal-order case}
A first $S$-minimum minimizes $S$ on the nonidentity norm-one units; a second successive $S$-minimum minimizes $S$ among units multiplicatively independent of the first. Norm normalization identifies these with successive minima in the full unit group modulo sign.

\begin{theorem}[Brunotte--Halter-Koch]\label{thm:bhk-index}
  Let $K$ be a totally real non-Galois cubic field and put $E=E^1(\OK)$. Choose $u\in E\setminus\{1\}$ with $S(u)=m_{\OK}$, and choose $v\in E\setminus\gp{u}$ minimizing $S$ there. Put $q=[E:\gp{u,v}]$. Then $q\le4$, and
  \[
    \begin{array}{c|c}
      q & \text{necessary relation for some }z\in E \\ \hline
      2 & uv=z^2                                    \\
      3 & u/v=z^3                                   \\
      4 & u/v=z^4.
    \end{array}
  \]
\end{theorem}
\begin{proof}
  By Section~\ref{sec:eventual-powers}, $u$ is primitive in $E$, so $v$ is a second successive $S$-minimum. These are the necessary implications of \cite[Satz~1(a)--(d)]{BHK1979}; that theorem gives additional size bounds which are not needed here.
\end{proof}

\begin{lemma}\label{lem:ratio-cube-two}
  If distinct norm-one units $u,v$ in a cubic field $K$ satisfy $u+v=\pm2$, then $v/u$ is not a cube of an integral norm-one unit.
\end{lemma}
\begin{proof}
  Suppose $v/u=z^3$, and put $t=\Tr_{K/\Q}(z)$, $s=\Tr_{K/\Q}(z^{-1})$. Newton sums give
  \[
    \Nm_{K/\Q}(1+z^3)=t^3+s^3-6ts+8=\Nm_{K/\Q}((u+v)/u)=\pm8.
  \]
  For the plus sign, set $c=t+s$. Then $c^3=3ts(c+2)$, with $c=-2$ impossible, and
  \[
    (t-s)^2=\frac{c^2(6-c)}{3(c+2)}.
  \]
  Nonnegativity requires the integer $c$ to belong to $\{-1,0,1,\ldots,6\}$. Integrality of $ts=c^3/[3(c+2)]$ leaves only $c=0,6$, giving $(t,s)=(0,0),(3,3)$. The corresponding polynomials are $X^3-1$ and $(X-1)^3$. A nonrational cubic element cannot satisfy either, while a rational norm-one $z$ is $1$ and would give $u=v$.

  For the minus sign, $c^3+16=3ts(c+2)$ and
  \[
    (t-s)^2=-\frac{(c-4)^2}{3}-\frac{32}{3(c+2)}.
  \]
  This is negative for $c\ge-1$. For $c\le-3$, put $h=-c-2\ge1$; three times the right side is $-(h+6)^2+32/h\le-17$. The value $c=-2$ contradicts the original equation. Thus neither sign is possible.
\end{proof}

\begin{lemma}\label{lem:maximal-affine-generation}
  Let $K$ be a totally real non-Galois cubic field and $E=E^1(\OK)$. If distinct $u,v\in\mathcal M_{\OK}$ are affinely equivalent, then $E=\gp{u,v}$.
\end{lemma}
\begin{proof}
  By Section~\ref{sec:eventual-powers}, distinct global minima are independent and primitive in $E$. Hence they satisfy the first and second minimum hypotheses of Theorem~\ref{thm:bhk-index}. Lemmas~\ref{lem:translate-triple} and~\ref{lem:translate-double} give the following two cases.

  Suppose first the affine class contains a triple $u,v,w$ of norm-one translates with $uvw=1$. The product $uv$ cannot be a square in $E$, since its inverse $w$ is a primitive minimum. Thus $q\ne2$. If $q=4$, the relation $u/v=z^4$ would imply $uv=(vz^2)^2$, again impossible. If $q=3$, write $u/v=z^3$ and $d=u-v\in\Z\setminus\{0\}$. Taking norms in $v(z^3-1)=d$ gives
  \[
    d^3=\Nm_{K/\Q}(z^3-1)=t^3-s^3,\qquad t=\Tr_{K/\Q}(z),\quad s=\Tr_{K/\Q}(z^{-1}).
  \]
  Fermat's theorem for exponent three \cite[Theorem~227, \S13.4]{HardyWright1960} forces $ts=0$. If $s=0$, then $d=t$ and $z^3-1=tz^2$, so $v=z^{-2},u=z$. If $t=0$, then $d=-s$ and $z^3-1=-sz$, so $v=z^{-1},u=z^2$. In either case a minimum is a proper square. Therefore $q=1$.

  Suppose next the affine class contains exactly two representatives, with $u+v=\pm2$. If $uv$ were a square, then $u,v$ would have the same sign at every embedding. For sum $2$, both would be positive, with $u_iv_i\le1$ and strict inequality at some embedding, contradicting $\Nm_{K/\Q}(uv)=1$. For sum $-2$, both would be totally negative, contradicting norm one in odd degree. This excludes $q=2,4$ as before. Lemma~\ref{lem:ratio-cube-two} excludes $q=3$. Again $E=\gp{u,v}$.

\end{proof}

\subsection{An index-two overorder with unchanged units}
An \emph{overorder} of an order $R$ is an order containing $R$ in the same number field. The next family shows that an affine minimum need not generate the maximal order. Its member with $k=1$ also supplies the maximal-order calculation needed for one of the nineteen small cases.
\begin{theorem}\label{thm:index-two-family}
  Let $k\in\Z$ with $k\ge1$, let $\beta$ satisfy $X^3-4kX^2-2X+4k=0$, and put
  \[
    e=\frac{\beta^2-4k\beta}{2},\qquad
    \Omega=\Z+\Z\beta+\Z e,\qquad G(k)=64k^4+13k^2+2.
  \]
  Then $\Omega$ is an order with $[\Omega:\Z[\beta]]=2$, $\Omega^\times=\Z[\beta]^\times$, and
  \[
    m_\Omega=16k^2+6,\qquad\mathcal M_\Omega=\{\beta-1,-\beta-1\}.
  \]
  If the odd part of $G(k)$ is squarefree, $\Omega$ is the maximal order. Under this hypothesis all maximal-order minima are affinely equivalent, although $\Z[\beta]$ is nonmaximal.
\end{theorem}
\begin{proof}
  The polynomial is the minus-midpoint polynomial with $T=4k$, so it is irreducible, totally real and non-Galois by~\eqref{eq:midpoint-parameters}. Its discriminant is $16G(k)$. Direct multiplication gives
  \begin{equation}\label{eq:index-two-ring}
    \beta^2=4k\beta+2e,\qquad
    \beta e=\beta-2k,\qquad
    e^2=e-k\beta+4k^2.
  \end{equation}
  Thus $\Omega$ is a finite free ring over $\Z$ and hence an order of algebraic integers. Its index above $R=\Z[\beta]$ is two, with $\disc(\Omega)=4G(k)$.

  The relations~\eqref{eq:index-two-ring} define two unital ring maps $\Omega\to\mathbb F_2$ by
  \[
    (\beta,e)\longmapsto(0,0),\qquad(\beta,e)\longmapsto(0,1).
  \]
  If $\alpha=n+x\beta+ye$ has odd $y$, its images are $n$ and $n+1$, one of which is zero. Therefore it cannot be a unit. Every $\Omega$-unit lies in $R$, and its inverse does too. The opposite inclusion is immediate, proving $\Omega^\times=R^\times$. Proposition~\ref{prop:midpoint-order} gives the stated minimum and complete pair.

  We prove that $\Omega$ is maximal at $2$, meaning that $2\nmid[\OK:\Omega]$. We use the following special case of Dedekind's index criterion \cite[Theorem~6.1.4(2)]{Cohen1993}: if a monic irreducible polynomial $P$ reduces modulo $2$ to $X(X+1)^2$ and $P(-1)/2$ is odd, its equation order is maximal at $2$. Indeed the only repeated irreducible factor is $X+1$; it does not divide $(P-X(X+1)^2)/2$ modulo $2$, as evaluation at $-1$ shows.

  If $k$ is odd, $e$ has polynomial
  \[
    P_e(X)=X^3-2X^2+(1-4k^2)X+2k^2.
  \]
  It reduces to $X(X+1)^2$ and satisfies $P_e(-1)=2(3k^2-2)$, of $2$-adic valuation one. The determinant of $1,e,e^2$ in the basis $1,\beta,e$ is $k$, so $[\Omega:\Z[e]]=k$ is odd. Hence $\Omega$ is maximal at $2$.

  If $k$ is even, use $\alpha=e+\beta$, whose polynomial is
  \[
    P_\alpha(X)=X^3-(4k+2)X^2+(-4k^2+10k-1)X
    +16k^3-14k^2+2k.
  \]
  Again its reduction is $X(X+1)^2$, and
  \[
    P_\alpha(-1)=-2(1+6k+5k^2-8k^3)
  \]
  has valuation one. Since $\alpha^2=(3k+2)\beta+3e+4k^2-4k$, the index $[\Omega:\Z[\alpha]]=3k-1$ is odd, proving maximality at $2$ here also. Both generators have degree three: rationality of either would impose a nonzero quadratic equation on $\beta$.

  Finally, if the odd part of $G(k)$ is squarefree, no odd prime square divides $\disc(\Omega)=4G(k)$, so no odd prime divides $[\OK:\Omega]$. The local argument excludes $2$. Thus $\Omega=\OK$.
\end{proof}

For $R=\Z[\beta]$, the conductor of $R$ in $\Omega$, the largest ideal of $\Omega$ contained in $R$, is $2\Z+\Z\beta+2\Z e$. The corresponding quotient of $\Omega$ is $\mathbb F_2\times\mathbb F_2$, with the image of $R$ diagonal; both residue unit groups contain only $(1,1)$. Thus the units are unchanged although $S(e)=12k^2+1$ is below their minimum. A minimizing unit need not generate the maximal order.

For $k=1$, $G(1)=79$, the field discriminant is $316$, and the norm-one minima are $\beta-1,-\beta-1$ at level $22$. This field already occurs in \cite[\S3]{Berwick1932}: if $\alpha^3-7\alpha-2=0$, then $\beta=(\alpha^2+\alpha)/2-1$. The polynomial family is a subfamily of \cite[Theorem~3]{Minemura1998}, with $\theta=-\beta-1$ and $t=4k+1$. The polynomial defining $\beta$ has discriminant divisible by $16$, so Minemura's squarefree-discriminant theorem does not give the proper-overorder conclusion proved here. The maximal-order conclusion requires the stated arithmetic hypothesis; no infinitude of squarefree odd parts of $G(k)$ is asserted.

\subsection{Proof for arbitrary rank-two subgroups}\label{sec:subgroup-minima}

For common minima above $41$, the individual-unit inequalities of Brunotte--Halter-Koch give the required index bound. The finitely many smaller affine forms require the full-group theorem.

\begin{lemma}\label{lem:relative-index}
  Let $\Gamma\subseteq E^1(\OK)$ have rank two in a totally real non-Galois cubic field. If distinct $u,v\in\mathcal M_\Gamma$ have common value $m>41$, then $r=[\Gamma:\gp{u,v}]\le4$. At indices $2,3,4$, respectively, there is a unit $z\in\Gamma$ with $uv=z^2$, $u/v=z^3$, or $u/v=z^4$.
\end{lemma}
\begin{proof}
  We adapt the argument of \cite[proof of Satz~1]{BHK1979} so that every auxiliary root belongs to $\Gamma$. By Section~\ref{sec:eventual-powers}, $u,v$ are primitive and independent. Extend $u$ to a basis $(u,h)$ of $\Gamma$ and write $v=u^ch^r$, with $r>0$. The exponent determinant gives $r=[\Gamma:\gp{u,v}]$, and primitivity of $v$ gives $\gcd(c,r)=1$. There is nothing to prove if $r=1$. If $r=2$, then $c$ is odd and $z=u^{(1+c)/2}h\in\Gamma$ satisfies $z^2=uv$. For $r\ge3$, choose $\delta\in\{1,-1\}$ and $1\le q\le r/2$ so that $q+\delta c$ is divisible by $r$. Here $q$ is an exponent, not a group index. Then $\gcd(q,r)=1$, $q\le(r-1)/2$, and
  \[
    \eta=u^{(q+\delta c)/r}h^\delta\in\Gamma,
    \qquad \eta^r=u^qv^\delta,\qquad S(\eta)\ge m.
  \]
  The $h$-exponent of $\eta$ is $\delta$, so $\eta\ne1$ and minimality gives the last inequality. The individual-unit estimates in \cite[Lemma~2]{BHK1979} give, for nonidentity integral norm-one units $x,y$,
  \begin{gather*}
    S(x^{-1})\le S(x)^2,\qquad S(xy)<3S(x)S(y),\\
    S(x^2)\le\frac{51}{32}S(x)^2\quad(S(x)>19),\qquad
    \frac{S(x)^3}{9}\le S(x^3)\le\frac{151}{49}S(x)^3,\\
    S(x)^n<\frac{3^{n+1}}2 S(x^n),\qquad
    S(x^n)<3\left(\frac32\right)^nS(x)^n
    \quad(S(x)>34,\ n\ge2).
  \end{gather*}
  The last upper bound also holds for $n=1$. Thus
  \begin{equation}\label{eq:relative-index-bound}
    m^r<\frac{3^{r+2}}2 S(u^q)S(v^\delta).
  \end{equation}

  Suppose first that $\delta=1$. If $r\ge8$, the general upper-power bound gives
  \[
    m<\frac32\left(\frac{3^{r+3}}2\right)^{2/(r-1)}
    \le\frac32\left(\frac{3^{11}}2\right)^{2/7}<39.
  \]
  For $r\le7$, one has $q\in\{1,2,3\}$. If $q=1,r=3$, the cubic lower bound gives $m^3<27m^2$, hence $m<27$. The remaining bounds from~\eqref{eq:relative-index-bound} are
  \[
    \begin{array}{c|c|c}
      q & \text{range} & \text{upper bound for }m                                    \\\hline
      1 & r\ge4        & (3^{r+2}/2)^{1/(r-2)}\le(3^6/2)^{1/2}<20                    \\
      2 & r\ge5        & (3^{r+2}\cdot51/64)^{1/(r-3)}\le(3^7\cdot51/64)^{1/2}<42    \\
      3 & r\ge7        & (3^{r+2}\cdot151/98)^{1/(r-4)}\le(3^9\cdot151/98)^{1/3}<32.
    \end{array}
  \]
  These estimates use, respectively, $S(u)=m$, the square bound, and the cube bound. Every displayed expression decreases with $r$: each is a constant times a fixed base greater than one raised to a decreasing positive exponent. Since $m$ is an integer greater than $41$, all positive-sign cases are impossible.

  Suppose next that $\delta=-1$. If $r\ge7$ and $q\le(r-3)/2$, the inverse and general power bounds give
  \[
    m<\frac12\,3^{(3r+3)/(r-1)}\le\frac{81}{2}<41.
  \]
  If $r\ge7$ and $q=(r-2)/2$, put $\xi=u\eta^{-2}\in\Gamma$. Then $\xi^r=(uv)^2$, so $\xi^{r/2}=uv$ since $\Gamma$ is torsion-free. The $h$-coordinate of $\xi$ is $2$, so $\xi\ne1$ and $S(\xi)\ge m$. The positive-sign $q=1$ estimate excludes this equation because $r/2\ge4$. That estimate uses only the root equation and a nonidentity root in $\Gamma$; the root exponent need not be the original subgroup index.

  If $r\ge5$ and $q=(r-1)/2$, the same choice $\xi=u\eta^{-2}$ satisfies $\xi^r=uv^2$. Again its $h$-coordinate is $2$, so $\xi\ne1$ and $S(\xi)\ge m$. The positive-sign $q=2$ estimate, with $u,v$ interchanged, gives $m<42$.

  The remaining possibilities have $q=1$ and $r\in\{3,4,5,6\}$. At $r\ge5$,~\eqref{eq:relative-index-bound} gives
  \[
    m<(3^{r+2}/2)^{1/(r-3)}\le(3^7/2)^{1/2}<34.
  \]
  Thus only $r=3,4$ remain, with $\eta^r=u/v$, as required.
\end{proof}

\pagebreak[3]
\begin{lemma}\label{lem:small-affine-subgroups}
  In the setting of Theorem~\ref{thm:subgroup-affine-minima}, if $m_\Gamma\le41$, then $\Gamma=\gp{u,v}$.
\end{lemma}
\begin{proof}
  We enumerate the affine normal forms at spread at most $41$, then determine whether a larger subgroup could retain their minimum. A pair of norm-one integral translates completes to a third translate $(uv)^{-1}\in\Gamma$; repetition would contradict primitivity. Lemmas~\ref{lem:translate-triple} and~\ref{lem:translate-double} therefore give either a product triple $\theta-b,\theta,\theta+a$ or the midpoint pair~\eqref{eq:midpoint-family}. In the triple case $a=1$ is excluded by
  \[
    S(\theta^{-1})=b^2-3b+3<b^2+b+1=S(\theta).
  \]
  If $b=1$, the non-Galois cases below level $42$ are the Ennola parameters $a=3,4,5$. If $a,b\ge2$, the bound $a^2+ab+b^2\le41$ forces $a,b\le5$ and gives
  \[
    (a,b)=(2,2),(2,3),(2,4),(3,2),(3,3),(3,4),(4,2),(4,3),(5,2).
  \]
  The only additional pair $(2,5)$ has square discriminant $4225$ and is excluded. The midpoint ranges are $\delta=1$, $T=4,5,6$, and $\delta=-1$, $T=2,3,4,5$. These are nineteen forms in total. Write $H$ for the subgroup generated by the displayed pair and $q_K=[E^1(\OK):H]$.

  Twelve forms have the following squarefree polynomial discriminants:
  \[
    \begin{array}{c|c|c@{\qquad}c|c|c}
      \text{family}  & \text{parameters} & \Delta     & \text{family}     & \text{parameters} & \Delta    \\\hline
      \text{product} & (2,2)             & 229        & \text{product}    & (5,2)             & 5521      \\
      \text{product} & (2,3)             & 761        & \text{Ennola}     & 3                 & 257       \\
      \text{product} & (2,4)             & 19\cdot103 & \text{Ennola}     & 4                 & 17\cdot41 \\
      \text{product} & (3,2)             & 5\cdot197  & \text{Ennola}     & 5                 & 1489      \\
      \text{product} & (3,4)             & 11\cdot619 & \text{midpoint }- & 3                 & 11\cdot43 \\
      \text{product} & (4,3)             & 11\cdot659 & \text{midpoint }- & 5                 & 2857
    \end{array}
  \]
  All factors displayed here are prime. Their equation orders are maximal, so the equation-order minima and Lemma~\ref{lem:maximal-affine-generation} give $q_K=1$. The seven remaining discriminants are
  \begin{equation}\label{eq:small-affine-table}
    \begin{array}{c|c|c}
      \text{family}     & \text{parameters} & \Delta             \\\hline
      \text{product}    & (3,3)             & 3^3\cdot107        \\
      \text{product}    & (4,2)             & 7^2\cdot53         \\
      \text{midpoint }+ & 4                 & 2^4\cdot37         \\
      \text{midpoint }+ & 5                 & 5^2\cdot73         \\
      \text{midpoint }+ & 6                 & 2^2\cdot3^4\cdot13 \\
      \text{midpoint }- & 2                 & 2^2\cdot37         \\
      \text{midpoint }- & 4                 & 2^4\cdot79
    \end{array}
  \end{equation}
  Four more equation orders are maximal. The following table lists a defining polynomial $g$ and every possible index prime $p$. For product $(4,2)$, we first translate by $X\mapsto X+4$.
  \[
    \begin{array}{c|c|c|c}
      \text{family}           & g(X)             & p   & g(0) \\\hline
      \text{product }(4,2)    & X^3+14X^2+56X+63 & 7   & 63   \\
      \text{midpoint }+,\ T=5 & X^3-5X^2+5       & 5   & 5    \\
      \text{midpoint }+,\ T=6 & X^3-6X^2+6       & 2,3 & 6    \\
      \text{midpoint }-,\ T=2 & X^3-2X^2-2X+2    & 2   & 2
    \end{array}
  \]
  In each test, $g(X)\equiv X^3\pmod p$ and $v_p(g(0))=1$, so Dedekind's index criterion excludes $p$ from the index \cite[Theorem~6.1.4(2)]{Cohen1993}. This includes the wild prime $3$ at $T=6$. The equation-order minima and Lemma~\ref{lem:maximal-affine-generation} therefore give $q_K=1$ in all four cases.

  For midpoint minus $T=4$, Theorem~\ref{thm:index-two-family} at $k=1$ gives the maximal order and its unchanged minimum pair, so Lemma~\ref{lem:maximal-affine-generation} gives $q_K=1$ here as well.

  For product $(3,3)$, let $\theta^3-9\theta-1=0$, work in $K=\Q(\theta)$, and set $z=(-\theta^2-\theta+5)/3$. Direct substitution gives
  \[
    z^3+z^2-4z-1=0,\qquad z^2=\theta+3,\qquad z^{-1}=\theta+z-1.
  \]
  The integral lattice with basis $1,\theta,z$ has trace Gram matrix
  \[
    \begin{pmatrix}3&0&-1\\0&18&-7\\-1&-7&9\end{pmatrix}
  \]
  of squarefree determinant $321$, so $\{1,\theta,z\}$ is an integral basis of $\OK$. Its spread form is
  \[
    Q(x,y)=27x^2-21xy+13y^2,\qquad
    52Q(x,y)=(26y-21x)^2+963x^2.
  \]
  The bound $Q\le27$ gives $|x|\le1$; direct substitution leaves, up to simultaneous sign, only the directions $z,z^{-1},\theta$, of levels $13,19,27$. The first two have norm polynomials
  \[
    \Nm_{K/\Q}(n+z)=n^3-n^2-4n+1,\qquad
    \Nm_{K/\Q}(n+z^{-1})=n^3-4n^2-n+1.
  \]
  Norm one requires $n=0$ in either case. For norm minus one, the rational-root test leaves only $n=\pm1,\pm2$, none of which works. Thus $z$ is a first full-group minimum, and $\theta$ is a second minimum outside $\gp z$. For independence, at the embedding $\theta>3$ one has $z<0$, so $\theta$ is not an odd power of $z$; the two negative conjugates of $\theta$ exclude even powers. Since $z$ is primitive, dependence would make $\theta$ an integral power of $z$.

  Apply Theorem~\ref{thm:bhk-index} to $(z,\theta)$. The negative sign of $z\theta$ and $z/\theta$ at $\theta>3$ excludes indices two and four. If $z/\theta=w^3$, the polynomial $X^3+16X^2+11X-1$ of $z/\theta$ gives, with $t=\Tr_{K/\Q}(w)$, $s=\Tr_{K/\Q}(w^{-1})$,
  \[
    t^3-3ts+3=-16,\qquad s^3-3ts+3=11,\qquad t^3-s^3=-27.
  \]
  Fermat's theorem for exponent three leaves $(t,s)=(-3,0),(0,3)$, neither satisfying the first equation. Therefore $E^1(\OK)=\gp{z,\theta}$. The identity $\theta-3=\theta^{-1}z^{-2}$ gives $H=\gp{\theta,z^2}$ and $q_K=2$.

  For midpoint plus $T=4$, put $z=(-\beta^2+2\beta+2)/2$. Then
  \[
    z^3+z^2-3z-1=0,\qquad z^2=\beta+1.
  \]
  The element $z+1$ has the already settled midpoint-minus $T=2$ polynomial, so $E^1(\OK)=\gp{z,h}$ with $h=-z-2$. If $v=1-\beta$, the identity $vh=z^{-1}$ gives
  \[
    H=\gp{\beta+1,1-\beta}=\gp{z^2,z^{-1}h^{-1}},
  \]
  whose exponent determinant has absolute value two. Hence $q_K=2$ here also.

  In each of the seventeen index-one forms, $H\subseteq\Gamma\subseteq E^1(\OK)=H$ gives equality. In either index-two form, strict containment $H\subsetneq\Gamma$ would force $\Gamma=E^1(\OK)$. But the displayed $z$ has spread $13<27$ in product $(3,3)$ and $10<16$ in midpoint plus $T=4$, contradicting $\Gamma$-minimality. Thus $\Gamma=H$ in every case.
\end{proof}

\begin{proof}[Proof of Theorem~\ref{thm:subgroup-affine-minima}]
  Section~\ref{sec:eventual-powers} gives primitivity and independence inside $\Gamma$. A pair of integral translates completes to $(uv)^{-1}\in\Gamma$; primitivity prevents repetition, so Lemma~\ref{lem:translate-triple} gives a product-one triple of $\Gamma$-minima. Otherwise the pair is reflected, and Lemma~\ref{lem:translate-double} gives $u+v=\pm2$.

  For $m_\Gamma\le41$, Lemma~\ref{lem:small-affine-subgroups} gives generation. For $m_\Gamma>41$, apply Lemma~\ref{lem:relative-index}. The exclusions in the proof of Lemma~\ref{lem:maximal-affine-generation} take place entirely inside $\Gamma$: in the triple case a square product contradicts primitivity of $(uv)^{-1}$, and a cube quotient gives $d^3=t^3-s^3$ and then a proper square among the minima. In the reflected case a square product is excluded by signs and norms, and Lemma~\ref{lem:ratio-cube-two} excludes a cube quotient. A fourth-root quotient would make the product a square. All roots supplied by the relative index lemma belong to $\Gamma$, so these arguments exclude indices two, three and four. Hence $\Gamma=\gp{u,v}$.

  It remains to determine all minima. Set $R=\Z[u]$. Affine equivalence and the integral polynomial expressions for unit inverses give $v,u^{-1},v^{-1}\in R$. Since $\Gamma=\gp{u,v}$, it follows that $\Gamma\subseteq E^1(R)$ and $\Z[\Gamma]=R$. The inverse calculation in the proof of Lemma~\ref{lem:small-affine-subgroups}, valid for every $b\ge1$, excludes $a=1$ in a product triple. If $b=1$, the case $a=2$ is cyclic and hence excluded. The remaining triples are covered by Theorem~\ref{thm:equation-order} and Proposition~\ref{prop:product-order}; the reflected pairs are covered by Proposition~\ref{prop:midpoint-order}. Each equation-order result gives minimum value $S(u)=m_\Gamma$ and a complete minimum set already contained in $\gp{u,v}$. Since $\Gamma\subseteq E^1(R)$, that set is exactly $\mathcal M_\Gamma$. Every member is affinely equivalent to $u$, so $R=\Z[w]$ for each minimum $w$. The product-one relation makes every pair in the triple generate the same group.
\end{proof}

\section{Orders, families, and the affine-pair hypothesis}\label{sec:applications}

\subsection{Unit-generated orders}
The subgroup theorem identifies the smallest order retaining the minimum and gives a generation criterion for the classical families.

\begin{corollary}\label{cor:subgroup-unit-order}
  Under the hypotheses of Theorem~\ref{thm:subgroup-affine-minima}, one has $\Gamma=E^1(\Z[\Gamma])$.
\end{corollary}
\begin{proof}
  Theorem~\ref{thm:subgroup-affine-minima} gives $R=\Z[\Gamma]=\Z[u]$ and $\Gamma=\gp{u,v}$. The equation-order calculations in its proof show that $u,v$ are also global minima in $E^1(R)$. The norm-one unit group of an order has rank two, so the same theorem applies to $E^1(R)$ and gives $E^1(R)=\gp{u,v}=\Gamma$.
\end{proof}

\begin{proof}[Proof of Corollary~\ref{cor:unit-generated-order}]
  Apply Theorem~\ref{thm:subgroup-affine-minima} to $E^1(\Omega)$, and use $\Omega^\times=\{\pm1\}E^1(\Omega)$. The additive assertion follows from the equation-order equality cases in Theorem~\ref{thm:equation-order} and Propositions~\ref{prop:product-order} and~\ref{prop:midpoint-order}. If $m_\Xi=m_\Omega$, a minimizing unit $w$ of $\Xi$ is an $\Omega$-minimum, so $R_U=\Z[w]\subseteq\Xi$. This inclusion puts every unit of $\Omega$ and its inverse in $\Xi$, giving equality of the unit groups. Equal groups have the same minimum and complete minimizing set.
\end{proof}

Thus $R_U$ is the unique smallest suborder retaining $\Omega$'s unit group or its minimum value. Its overorders inside $\Omega$ have the same regulator and complete minimum set, although their additive minima can differ. For arbitrary subgroups, equality of minimum values need not imply equality of the subgroups or of their complete minimum sets. If $u,v$ form a basis of $E^1(\Omega)$ consisting of global minima as above, then $\Gamma=\gp{u,v^2}$ has index two and $m_\Gamma=m_\Omega$, whereas $\mathcal M_\Gamma=\{u\}$. The other minima are $v$ and, in the triple case, $(uv)^{-1}$, both with odd $v$-exponent.

\subsection{The classical families}
\begin{theorem}\label{thm:family-generation}
  Let $R$ be one of the orders in Theorem~\ref{thm:equation-order} or Propositions~\ref{prop:product-order} and~\ref{prop:midpoint-order}, and suppose its field of fractions $K$ is non-Galois. Let $M_0$ be the displayed minimum value for $R$, and let $H$ be generated by two distinct displayed norm-one minima. For every subgroup $H\subseteq\Gamma\subseteq E^1(\OK)$, one has $m_\Gamma=M_0$ if and only if $\Gamma=H$. In this case the complete minimum set is the displayed pair or triple; otherwise $m_\Gamma<M_0$. In particular, $E^1(R)=H$ without any maximality hypothesis on $R$.
\end{theorem}
\begin{proof}
  The equation-order calculations give $m_H=M_0$ and the displayed equality set. The displayed units are global $E^1(R)$-minima and hence independent by Section~\ref{sec:eventual-powers}, so every intermediate $\Gamma$ has rank two. They give $m_\Gamma\le M_0$. Equality makes them global $\Gamma$-minima, so Theorem~\ref{thm:subgroup-affine-minima} gives $\Gamma=H$. The converse follows from the equation-order calculation. Taking $\Gamma=E^1(R)$ proves the last assertion.
\end{proof}

For the Ennola family, put $E_a=E^1(\mathcal O_a)$ and $A_a=a^2+a+1$. Theorem~\ref{thm:family-generation} says that every subgroup $H_a\subseteq\Gamma\subseteq E_a$ satisfies $m_\Gamma=A_a$ if and only if $\Gamma=H_a$. The complete minimum set is then $\{\theta-1,\theta,\theta+a\}$; otherwise $m_\Gamma<A_a$. This applies to every order $R_a\subseteq\Omega\subseteq\mathcal O_a$ by taking $\Gamma=E^1(\Omega)$.

The remaining question for Ennola fields is whether all maximal-order minima are affinely equivalent when their spread is below $A_a$. Ennola's conjecture $E_a=H_a$ would imply that they are affinely equivalent by the preceding theorem. Choi and Kim show that, for each fixed prime $p\ge3$, only finitely many $a\ge3$ satisfy $p\mid[E_a:H_a]$ \cite[Theorem~1]{ChoiKim2025}. The exceptional set may depend on $p$, so this does not establish eventual index one. The subgroup theorem does not depend on these open questions.

\subsection{Minima in distinct affine classes}\label{sec:counterexample}

To show why the affine-pair hypothesis is needed, we give a maximal order with minima in two different affine classes. The proof bounds the possible additive classes by a positive definite binary quadratic form and checks the norm equations in every remaining class.

\begin{proposition}\label{prop:general-counterexample}
  Let $\theta$ satisfy $X^3-14X^2-14X+14=0$ and put $K=\Q(\theta)$. Then $\OK=\Z[\theta]$, $D_K=247156$, and $m_{\OK}=952$. There are exactly two minimizing affine classes, represented by
  \[
    u=-2\theta-3,\qquad v=\theta^2-15\theta+9.
  \]
  Both displayed units have norm one and are not affinely equivalent.
\end{proposition}
\begin{proof}
  The polynomial is Eisenstein at $7$, with discriminant $2^2\cdot7^2\cdot13\cdot97$. It is also Eisenstein at $2$. Both primes define tame totally ramified cubic extensions, with discriminant exponent two. Since only $2,7$ could divide the equation-order index, that index is one. The positive nonsquare discriminant makes the field totally real and non-Galois.

  Set $\eta=\theta^2-14\theta$. In the quotient basis $\theta,\eta$,
  \[
    Q(x,y)=238x^2+70xy+784y^2,\qquad
    952Q(x,y)=(476x+70y)^2+741468y^2.
  \]
  Thus $Q\le952$ implies $|y|\le1$. Among nonzero quotient classes, the complete list modulo simultaneous sign is
  \[
    (1,0),\quad(0,1),\quad(2,0),\quad(-1,1),
  \]
  at levels $238,784,952,952$.

  The shorter classes have norm polynomials
  \[
    \Nm_{K/\Q}(n+\theta)=n^3+14n^2-14n-14,\qquad
    \Nm_{K/\Q}(n+\eta)=n^3+28n^2-2548.
  \]
  The first takes only residues $0,2,4,6,7,8,9$ modulo $11$, excluding both unit norms. For the second, the complete residue sets are
  \[
    \begin{array}{c|l}
      11 & 0,3,4,7,8,9,10,                 \\
      19 & 1,2,3,4,5,6,7,8,11,12,14,16,17.
    \end{array}
  \]
  These exclude $1$ modulo $11$ and $-1$ modulo $19$, respectively.

  At the remaining two classes,
  \[
    \begin{aligned}
      \Nm_{K/\Q}(n+2\theta)     & =n^3+28n^2-56n-112,  \\
      \Nm_{K/\Q}(n+\eta-\theta) & =n^3+14n^2-252n+406.
    \end{aligned}
  \]
  The first has value $-1$ at $n=3$, and the second has value $1$ at $n=9$, giving $u,v$. Their $\theta^2$ coefficients differ by one, whereas that of every $\pm u+n$ is zero. Hence they are inequivalent. The exhaustive four-class list proves global minimality and the precise number of minimizing classes.
\end{proof}

By Theorem~\ref{thm:subgroup-affine-minima}, each of these two affine classes contains only one minimum: a second minimum in either class would force every minimum into that class. Thus $\mathcal M_{\OK}=\{u,v\}$.

The LMFDB records this field and lists $u,v$ as fundamental units \cite{LMFDB247156}. With $\alpha=5-\theta$, the database's defining polynomial is $X^3-X^2-79X+281$, and its listed units are $2\alpha-13$ and $\alpha^2+5\alpha-41$. The quadratic-form calculation and Theorem~\ref{thm:subgroup-affine-minima} together show that these are the complete global minimum set for $S$.

There is also a distinction between additive and multiplicative primitivity in this example. As a global minimum, $u$ is multiplicatively primitive by Section~\ref{sec:eventual-powers}, but its coordinates in $\OK/\Z$ have greatest common divisor two. Thus a primitive unit need not represent a primitive additive class.

\renewcommand{\bibliofont}{\fontsize{9}{10.5}\selectfont}
\bibliographystyle{alpha}
\bibliography{references}

\end{document}